\documentclass[onecolumn]{autart}
\usepackage{savesym}
\savesymbol{AND}
\usepackage{graphicx}
\usepackage{amssymb}
\usepackage{amsmath}
\usepackage{epstopdf}
\usepackage{dsfont}
\usepackage{xcolor,colortbl}
\usepackage{algorithm}
\usepackage{algorithmic}
\usepackage{todonotes}
\usepackage{hyperref}
\usepackage{booktabs}
\usepackage{caption}
\usepackage{mathtools}
\usepackage{comment}
\newcommand {\be}{\begin{equation}}
\newcommand {\ee}{\end{equation}}

\newcommand {\beW}{\begin{equation*}}
\newcommand {\eeW}{\end{equation*}}

\newcommand {\bsp}{\begin{split}}
\newcommand {\esp}{\end{split}}

\newcommand {\bea}{\begin{eqnarray}}
\newcommand {\eea}{\end{eqnarray}}

\newcommand {\beaW}{\begin{eqnarray*}}
\newcommand {\eeaW}{\end{eqnarray*}}

\newcommand {\bM}{\begin{bmatrix}}
\newcommand {\eM}{\end{bmatrix}}

\newcommand {\bi}{\begin{itemize}}
\newcommand {\ei}{\end{itemize}}

\newcommand {\ben}{\begin{enumerate}}
\newcommand {\een}{\end{enumerate}}

\newcommand {\bal}{\begin{align}}
\newcommand {\eal}{\end{align}}

\newcommand {\balW}{\begin{align*}}
\newcommand {\ealW}{\end{align*}}

\newcommand {\bcm}{\begin{columns}}
\newcommand {\ecm}{\end{columns}}
\newcommand {\bc}{\begin{column}}
\newcommand {\ec}{\end{column}}

\newcommand{\m}{\mathcal}
\newcommand{\mb}{\mathbf}

\newtheorem{assumption}{Assumption}

\DeclareMathOperator{\dx}{dx}

\DeclareMathOperator{\dz}{dz}

\DeclareMathOperator{\Xs}{\mathcal{X}}
\DeclareMathOperator{\F}{\mathcal{F}}

\DeclareMathOperator{\XUs}{\mathcal{X}\times\mathcal{U}}

\DeclareMathOperator{\Us}{\mathcal{U}}
\DeclareMathOperator{\Prb}{\mathbb{P}}

\DeclareMathOperator{\expect}{\mathbb{E}}

\DeclareMathOperator{\Ne}{\mathbb{N}}

\DeclareMathOperator{\st}{subject\ to}
\DeclareMathOperator{\trace}{tr}
\DeclareMathOperator{\sthat}{such\ that}

\DeclarePairedDelimiter\parens{\lparen}{\rparen}
\DeclarePairedDelimiter\bracks{\{}{\}}

\begin{document}
%%%%%%%%%%%%%%%%%%%%%%%%%%%%%%%%%%%%%%%%%%%%%%%%%%%%%%%%%%%%%%%%%%%%%%%%%%%%%%%%%%%%%%%%%%%%%%%%%%%%%%%%%%%%%%%%%%%%%%%%%%%%%%%%%%%%%%%%%%%%%%%%
%%BEGIN PAPER HERE
%%%%%%%%%%%%%%%%%%%%%%%%%%%%%%%%%%%%%%%%%%%%%%%%%%%%%%%%%%%%%%%%%%%%%%%%%%%%%%%%%%%%%%%%%%%%%%%%%%%%%%%%%%%%%%%%%%%%%%%%%%%%%%%%%%%%%%%%%%%%%%%%
\begin{frontmatter}
\title{Upper bounds for the reach-avoid probability via robust optimization}

\thanks[footnoteinfo]{The work of N. Kariotoglou was supported by the Swiss National Science Foundation under grant number $200021\_137876$.}

\author[First]{Nikolaos Kariotoglou}\ead{karioto@control.ee.ethz.ch}, 
\author[First]{Maryam Kamgarpour}\ead{mkamgar@control.ee.ethz.ch}, 
\author[First]{Tyler H. Summers}\ead{tsummers@control.ee.ethz.ch}, and
\author[First]{John Lygeros}\ead{lygeros@control.ee.ethz.ch}
%\author{Nikolaos Kariotoglou, Maryam Kamgarpour \\Tyler Summers and John Lygeros
%\thanks{The authors are with the Automatic Control Laboratory, Department of Information Technology and Electrical Engineering, ETH Z\"urich, Z\"urich 8092, Switzerland (e-mail: karioto@control.ee.ethz.ch; mkamgar@control.ee.ethz.ch; tsummers@control.ee.ethz.ch; lygeros@control.ee.ethz.ch)}\thanks{The work of N. Kariotoglou was supported by the Swiss National Science Foundation under grant number $200021\_137876$.}\vspace{-1\baselineskip}
%}
%\date{\today}     
\address[First]{Automatic Control Laboratory, Department of Information Technology and Electrical Engineering, ETH Z\"urich, Z\"urich 8092, Switzerland}

\begin{keyword} % Five to ten keywords,
Reachability, approximate dynamic programming, value function bounds, Markov decision processes, stochastic control.
\end{keyword} % keyword list or with the

%%%%%%%%%%%%%%%%%%%%%%%%%%%%%%%%%%%%%%%%%%%%%%%%%%%%%%%%%%%%%%%%%%%%%%%%%%%%%%%%%%%%%%%%%%%%%%%%%%%%%%%%%%%%%%%%%%%%%%%%%%%%%%%%%%%%%%%%%%%%%%%%%%%
\begin{abstract}
We consider finite horizon reach-avoid problems for discrete time stochastic systems. Our goal is to construct upper bound functions for the reach-avoid probability by means of tractable convex optimization problems. We achieve this by restricting attention to the span of Gaussian radial basis functions and imposing structural assumptions on the transition kernel of the stochastic processes as well as the target and safe sets of the reach-avoid problem. In particular, we require the kernel to be written as a Gaussian mixture density with each mean of the distribution being affine in the current state and input and the target and safe sets to be written as intersections of quadratic inequalities. Taking advantage of these structural assumptions, we formulate a recursion of semidefinite programs where each step provides an upper bound to the value function of the reach-avoid problem. The upper bounds provide a performance metric to which any suboptimal control policy can be compared, and can themselves be used to construct suboptimal control policies. We illustrate via numerical examples that even if the resulting bounds are conservative, the associated control policies achieve higher reach-avoid probabilities than heuristic controllers for problems of large state-input space dimensions (more than 20). The results presented in this paper, far exceed the limits of current approximation methods for reach-avoid problems in the specific class of stochastic systems considered. 
\end{abstract}
\end{frontmatter}

%\begin{IEEEkeywords}
%Reachability, approximate dynamic programming, stochastic control, Markov decision processes.
%\end{IEEEkeywords}
%%%%%%%%% General Comments %%%%%%%%%%%%%%%%%
%%%%%%%%%%%%%%%%%%%%%%%%%%%%%%%%%%%%%%%%%%%%%%%%%%%%%%%%%%%%%%%%%%%%%%%%%%%%%
%\tyler{robust $\rightarrow$ infinite set}
%\tyler{Presentation of examples}
%
%\nikos{Something is wrong with scenario references. I dont find the theorems 3.5 and 3.3 where they should be.}
%\nikos{It would be interesting to show that this method can outperform gridding. Either by comparing with a coarsely gridded $\Us$ space or even one where both $\Xs,\Us$ are finite and show via Monte Carlo simulations that the suggested RBF approach performs better. Perhaps this is too much if we include all other examples - let me know what you think.}
%\nikos{Should we add a discussion with drawn conclusions from the numerical examples? If so, what would be the main messages? Are we going to comment on complexity?}
%\nikos{Rewrite algorithms with input data etc.}
%%%%%%%%%%%%%%%%%%%%%%%%%%%%%%%%%%%%%%%%%%%%%%%%%%%%%%%%%%%%%%%%%%%%%%%%%%%%%

%%%%%%%%%%%%%%%%%%%%%%% INTRODUCTION %%%%%%%%%%%%%%%%%%%%%%%%%%%%%%%%%%%%%%
\section{Introduction}
\label{sec:intro}
We deal with Markov decision processes of specific structure and focus on the stochastic reach-avoid problem where one maximizes the probability of reaching a target set while staying in a safe set, under given control constraints. This type of optimal control problem can be solved by a dynamic programming recursion \cite{abate2008probabilistic,summers2010verification} but its practical applications are limited due to the infamous curse of dimensionality that affects state-of-the-art space-gridding techniques. Despite their limitations, space-gridding methods are theoretically attractive since they provide straightforward ways to estimate approximation errors with respect to the optimal solution function, under general Lipschitz continuity assumptions on the associated system dynamics \cite{abate2007computational,prandini2006stochastic,kushner2001numerical}. In an attempt to extend the possible applications of the reach-avoid problem, we follow a less standard approach \cite{de2003linear} than dynamic programming to characterize the optimal value function of the stochastic reach-avoid problem.  We then construct an upper bound for the optimal value function and evaluate in practice the performance of the associated approximate control policy .

The main tools used in this work stem from convex optimization and methods used to reformulate robust inequalities \cite{wang2014approximate,bertsimas2006tractable,boyd2004convex} into semidefinite constraints. The authors in \cite{wang2014approximate} used a similar framework and showed its success in terms of control performance in a different type of optimal control problem. The framework considered here is fundamentally different due to the different type of basis approximation functions used and the way that reach-avoid problems handle state and control constraints naturally via the definition of their cost function. Another approach that addresses the approximation of reach-avoid problems based on the method suggested in \cite{de2003linear} is presented in \cite{kariotoglou2013approximate,kariotoglou2014adp}. Compared to the results in \cite{kariotoglou2014adp}, the technique developed here provides deterministic upper bounds to the reach-avoid value function (as opposed to probabilistic ones), scales to significantly higher dimensional problems, but requires stronger restrictions on the kernel of the stochastic processes. 

Our results are based on restricting the decision space of an infinite dimensional linear program \cite{de2003linear,kariotoglou2014adp} to the space spanned by a finite collections of Gaussian radial basis functions (RBFs) and then searching for means, variances and weights of a sum of such functions to upper bound the reach-avoid probability. We show that for some process kernels, such as the kernel of a linear system with additive Gaussian mixture noise, this restriction allows us to analytically compute the expectation of the stochastic process over the basis functions and reduce the problem of upper bounding the optimal reach-avoid value function to determining the parameters of a Gaussian RBF sum that upper bounds another such sum. We derive sufficient conditions in terms of linear matrix inequalities (LMIs) to characterize dominance between Gaussian RBF sums and then use them to formulate a semidefinite program (SDP) with which one can construct an upper bound to the optimal value function at each step of the reach-avoid recursion. We show via numerical examples that approximate control policies defined using these bounds scale to problems of significantly higher dimensions than what is currently possible with state gridding techniques.

The rest of the paper is organized as follows: In Section \ref{sec:stochreach_DP} we recall the basic definition of the reach-avoid problem for Markov decision processes and formulate a robust optimization problem whose solution is an upper bound to the reach-avoid value function. We focus most of Section \ref{sec:rob_reform} to reformulate the constraints of the robust optimization problem into sufficient conditions that can be checked efficiently with existing convex optimization software. Using these conditions we provide an algorithm to construct upper bounds for the value function of every step of a reach-avoid recursion. We then use the constructed bounds to define approximate control policies that maximize the reach-avoid probability. In Section \ref{sec:examples} we provide numerical examples to estimate the suboptimality gap between the constructed bounds and the optimal value function (computed via gridding) and illustrate that the approximate control policies can outperform heuristic controllers even when the bounds are not tight.
%%%%%%%%%%%%%%%%%%%%%%%%%%%%%%%%%%%%%%%%%%%%%%%%%%%%%%%%%%%%%%%&
%%%%%%%%%%%%%%%%%%%%%%%%%%%%%%%%%%%%%%%%%%%%%%%%%%%%%%%%%%%%%%%&
%%%%%%%%%%%%%%%%%%%%%%% Reachabuility %%%%%%%%%%%%%%%%%%%%%%%%%%
%%%%%%%%%%%%%%%%%%%%%%%%%%%%%%%%%%%%%%%%%%%%%%%%%%%%%%%%%%%%%%%&
\section{Stochastic Reach-Avoid Problem}
\label{sec:stochreach_DP}
%%%%%%%%%%%%%%%%%%%%%%
\subsection{Dynamic programming approach}
\label{sec:stochreach}
%%%%%%%%%%%%%%%%%%%%%%

% dynamics
Let $\Xs=\mathbb{R}^n$ denote a continuous state space and $\Us\subset \mathbb{R}^m$ a continuous and compact control space. We consider a discrete-time controlled stochastic process $x_{t+1}\sim Q(dx|x_t,u_t)$, $(x_t,u_t)\in \Xs \times\Us $ with a transition kernel $Q:\mathcal{B}(\Xs)\times\Xs\times\Us\rightarrow [0,1]$ where $\mathcal{B}(\mathcal{X})$ denotes the Borel $\sigma$-algebra of $\mathcal{X}$. Given a state control pair $(x_t,u_t)\in \mathcal{X}\times\mathcal{U}$, $Q(A|x_t,u_t)$ measures the probability of $x_{t+1}$ being in a set $A\in\mathcal{B}(\Xs)$. The transition kernel $Q$ is a Borel-measurable stochastic kernel, that is, $Q(A|\cdot)$ is a Borel-measurable function on $\m{X} \times \m{U}$ for each $ A \in \m{B}(\m{X})$ and $Q(\cdot|x,u)$ is a probability measure on $\m{X}$ for each $(x,u)$. For the rest of the paper all measurability conditions refer to Borel measurability. 
%We allow the state space $\mathcal{X}$ to be any subset of $\mathbb{R}^n$ and assume that the control space $\mathcal{U}\subseteq \mathbb{R}^m$  is compact; extensions to hybrid state and input spaces where some states or inputs are finite valued are possible \cite{abate2008probabilistic}. 
% reach avoid objective
We consider a safe set $K^{\prime}\in\mathcal{B}(\mathcal{X})$ and a target set $K\subseteq K^{\prime}$. We define an admissible $T$-step control policy to be a sequence of measurable functions $\pi = \{\pi_0,\dots,\pi_{T-1}\}$ where $\pi_i:\mathcal{X}\rightarrow \mathcal{U}$ for each $i\in\{0,\dots,T-1\}$. The reach-avoid problem over a finite time horizon $T$ is to find an admissible $T$-step control policy that maximizes the probability of $x_t$ reaching the set $K$ at some time $t_K\leq T$ while staying in $K^{\prime}$ for all $t\leq t_K$. For any initial state $x_0$ we denote the reach-avoid probability associated with a given $\pi$ as: $r_{x_0}^\pi(K,K^{\prime}) = \Prb^{\pi}_{x_0}\{\exists j \in [0,T] : x_j \in K \wedge \forall i \in [0,j-1],\ x_i \in K^{\prime}\setminus K \}$ and operate under the
convention that $[0,-1]=\emptyset$, which implies that the requirement on $i$ is automatically satisfied
when $x_0\in K$. 

In \cite{summers2010verification}, $r_{x_0}^\pi(K,K^{\prime})$ is shown to be equivalent to the following sum multiplicative cost function:

\begin{align}
r_{x_0}^\pi(K,K^{\prime}) = \mathbb{E}_{x_0}^\pi\left [\sum^{T}_{j=0}{\left( \prod^{j-1}_{i=0}{\mathds{1}_{K^{\prime}\setminus K}(x_i)}       \right)} \mathds{1}_K(x_j)        \right]
\label{eq:raprob}
\end{align} 
where $\prod_{i=k}^j(\cdot)=1$ if $k > j$. The function $\mathds{1}_A(x)$ denotes the indicator function of a set $A\in \mathcal{B}(\Xs)$ with $\mathds{1}_A(x)=1$ if $x\in A$ and $\mathds{1}_A(x)=0$ otherwise. The sets $K$ and $K^\prime$ can be time dependent or even stochastic \cite{summers2013stochastic} but for simplicity we assume here that they are constant. We denote the difference between the safe and target sets by $\bar{\Xs}:=K^\prime \setminus K$ to simplify the presentation of our results.

The solution to the reach-avoid problem is given by a dynamic recursion \cite{summers2010verification}. Define $V_k^*:\mathcal{X}\rightarrow [0,1]$ for $k=T-1,\dots,0$ by :
\begin{align}
\label{eq:DP}
\begin{split}
&V_k^*(x)=\sup_{u\in \mathcal{U}}\left\{\mathds{1}_K(x)+\mathds{1}_{\bar{\Xs}}(x)\int_{\mathcal{X}}{V^*_{k+1}(y)Q(dy|x,u)}\right\}\\
&V_{T}^*(x)=\mathds{1}_K(x).
\end{split}
\end{align}
The value of the above recursion at $k=0$ and for any initial state $x_0$ is the supremum of \eqref{eq:raprob} over all admissible policies, i.e. $V_0^*(x_0)=\sup_{\pi} r^\pi_{x_0}(K,K^\prime)$ \cite{summers2010verification}. In \cite{kariotoglou2014adp} the authors proved measurability of the value functions (and hence attainability of the supremum) in \eqref{eq:DP} under  the assumption that $Q$ is continuous with respect to  $u\in \Us$ for any given $x\in\Xs$. For the rest of this work we operate under the same kernel continuity assumption:

\begin{assumption}[Kernel continuity]
\label{ass:kernel_continuity}
For every $x \in \m{X}, A \in \m{B}(X)$ the mapping $u \mapsto Q(A|x,u)$ is continuous.
\end{assumption}

We define the optimal reach-avoid feedback policy at each state $k$ by
\begin{align}
\pi^*_k(x)=\arg\max_{u\in \mathcal{U}}\left\{\mathds{1}_K(x)+\mathds{1}_{\bar{\Xs}}(x)\int_{\mathcal{X}}{V^*_{k+1}(y)Q(dy|x,u)}\right\}.
\label{eq:optpol}
\end{align}
The only established way to approximate the solution of \eqref{eq:DP} is by gridding $\XUs$ and computing it backwards, starting from the known value function $V_{T}^*(x)=\mathds{1}_K(x)$. In this way, the value of $V_0^*(x)$ is approximated on the grid points of $\mathcal{X}$ while an approximate control policy at each step $k$ is computed by taking the maximum of \eqref{eq:optpol} over the grid points of $\mathcal{U}$. The advantages of this approach are that it is straightforward to estimate the approximation accuracy as a function of the grid size (under suitable regularity assumptions \cite{abate2007computational,prandini2006stochastic}) and that the approximate feedback control policy can be stored as a look-up table rather than computed online at each state. The disadvantage is that it quickly becomes intractable as the dimensions of the state and control spaces increase. In what follows we formulate an infinite dimensional linear program equivalent to the reach-avoid recursion and use tools from robust optimization to construct upper bound functions.
%In the following section we approximate the reach-avoid value function defined in (\ref{eq:DP}) in a different way by reformulating the dynamic recursion in an equivalent infinite dimensional LP. 
%%%%%%%%%%%%%%%%%%%%%%%%%%%%%%%%%%%%%%%%%%%%%%%%%%%%%%%%%%%%%%%%%%%%%%%%%%%%%
%%%%%%%%%%%%%%%%%%%%%%%%%% Linear Programming approach %%%%%%%%%%%%%%%%%%%%%&
%%%%%%%%%%%%%%%%%%%%%%%%%%%%%%%%%%%%%%%%%%%%%%%%%%%%%%%%%%%%%%%%%%%%%%%%%%%%%

\subsection{Robust optimization problem}
Inspired by the linear programming approach to dynamic programming \cite{de2003linear}, the authors in \cite{kariotoglou2014adp} formulated very similar version to the following infinite dimensional linear program over the space of measurable functions $\F=\{f:\Xs\rightarrow\mathbb{R},\ \text{$f$ measurable}\}$ to characterize the reach-avoid value function at time $k$:
\begin{align}
\label{eq:inflp_chap_rob}
\begin{split}
\min_{V(\cdot)\in \F} \quad &\int_{{\Xs}} V(x) \nu(\dx) \\
	\st \quad & V(x)\geq \mathds{1}_K(x)+\mathds{1}_{\bar{\Xs}}(x)\int_{\m{X}}V^*_{k+1}(z)Q(\dz|x,u), \quad \forall (x,u)\in \bar{\mathcal{X}}\times\mathcal{U}
\end{split}
\end{align}
where we use $\min$ instead of $\inf$ since, as shown in \cite{kariotoglou2014adp},  Assumption \ref{ass:kernel_continuity} ensures attainability of the optimizer.
The constraints in \eqref{eq:inflp_chap_rob} are equivalent to the pair of constraints 
\begin{align*}
V(x)&\geq 1,\ \forall x\in K \\
V(x)&\geq \int_{\m{X}}V^*_{k+1}(z)Q(\dz|x,u), \ \forall (x,u)\in \bar{\mathcal{X}}\times\mathcal{U}.
\end{align*}
The measure $\nu$ is interpreted as a state relevance measure \cite{de2003linear} that can be chosen to bias the value of the cost function in different places of the state space. In order to search for a feasible point in \eqref{eq:inflp_chap_rob}, we restrict the infinite dimensional function space $\mathcal{F}$ to a finite dimensional subspace $\tilde{\mathcal{F}}$. Assume that $V^*_{k+1}$ is known and consider the following semi-infinite linear program defined over functions $V\in \tilde{\F}\subset \F$ 
\begin{align}
	\label{eq:semiinfLP_chap_rob}
	\begin{split}
	\min_{V\in \tilde{\F}}\quad &\int_{{\Xs}} V(x)\nu(\dx)\\
	\st \quad & V(x)\geq \int_{\m{X}}V^*_{k+1}(z)Q(\dz|x,u), \quad \forall (x,u)\in \bar{\mathcal{X}}\times\mathcal{U}\\
	 \quad& V(x) \geq 1,\quad \forall x\in K.
	\end{split}
\end{align}
To guarantee that the minimum is attained in \eqref{eq:semiinfLP_chap_rob}, additional assumptions may be required on $\tilde{\F}$. For instance, $\tilde{\F}$ would have to be bounded to exclude situations where the optimization program has an unbounded cost function, e.g. $-\infty$ on the complement of $\bar{\Xs}$ which is unconstrained. An alternative approach is to explicitly add constraints on the value of $V$ on the complement of $\bar{\Xs}$, such as $V(x)\geq 0, \forall x\in\Xs\setminus\bar{\Xs}$ since $V$ represents a reach-avoid probability. The rest of the section is devoted to reformulating the constraints of \eqref{eq:semiinfLP_chap_rob} into conditions that can be checked efficiently using convex optimization software.

\section{Robust constraint reformulation}
\label{sec:rob_reform}
Checking the constraints in \eqref{eq:semiinfLP_chap_rob} numerically is, in general, intractable due to the infinite cardinality of the constraint set. Depending on the structure of $K$, $\bar{\Xs}\times \Us$, the function  $\int_{\m{X}}V^*_{k+1}(z)Q(\dz|x,u)$ and the representation of $V$ in $\tilde{\F}$, it may be possible to reformulate the constraints. In this section we impose specific structure to these elements and derive sufficient conditions to construct feasible solutions for \eqref{eq:semiinfLP_chap_rob}.
%\begin{assumption}
%We impose the following assumptions on the structure of the reach-avoid reachability problem:
%\begin{enumerate}
%\item The target and safe sets $K$ and $K^\prime$ can be written as finite unions of disjoint hyper-rectangle sets, i.e. 
%\begin{align*}
%K=\bigcup_{p=1}^P K_p=\bigcup_{p=1}^P \parens*{\bigtimes_{l=1}^n [a^p_l,b^p_l]},\quad  K^\prime=\bigcup_{m=1}^M K^\prime_m=\bigcup_{m=1}^M \parens*{\bigtimes_{l=1}^n [c^m_l,d^m_l]}
%\end{align*}
%for some finite $P,M \in\mathbb{N}_+$ with $n=\dim\parens*{\Xs}$ and $a^p,b^p,c^m,d^m\in \mathbb{R}^n$ for every $p$ and $m$.
%\item The state-relevance measure $\nu$ can be written as a product measure, i.e. $\nu(\dx)=\prod_{l=1}^n \nu_l(\dx_l)$.
%\end{enumerate}
%\label{assum:setup_impl}
%\end{assumption}
%

\subsection{Gaussian radial basis functions}

Let the restricted decision space  be $\tilde{\F}=\F^{M}$ where $\F^M$ denotes the span of a set $\{\phi(x,\mu_i,\Sigma_i)\}_{i=1}^{M}$ of $M\in\mathbb{N}_+$ Gaussian radial basis functions (RBFs) each one defined as,
\begin{align}
\label{eq:basisf_sdp}
\begin{split}
\phi(x,\mu_i,\Sigma_i):=\frac{1}{\sqrt{(2\pi)^n|\Sigma_i|}}\exp\left[-\frac{1}{2}(x-\mu_i)^{\top} {\Sigma_i}^{-1}(x-\mu_i)\right]
\end{split}
\end{align}
where the operator $|\cdot|$ denotes the determinant of the argument matrix.
Note that each function $\phi(x,\mu_i,\Sigma_i)$ corresponds to the density function of a multivariate Gaussian random variable $\mathcal{N}(\mu_{i},\Sigma_{i})$ with $x,\mu_{i}\in \mathbb{R}^n,\Sigma_{i}\in  \mb{S}_+^n$ where $\mb{S}_+^n$ denotes the space of (symmetric) positive definite matrices. We use the standard notation $X\succcurlyeq 0$ to denote that a matrix $X \in \mathbb{R}^{n\times n}$ is in $\mb{S}_+^n$ and hence we have that $\Sigma_i\succcurlyeq 0$. Every element of $\psi\in\F^{M}$ can be written as a Gaussian RBF sum defined as:
\begin{defn}[Gaussian RBF sum] A Gaussian RBF sum $\psi:\Xs\rightarrow \Re$ is a linear combination of functions in the form of \eqref{eq:basisf_sdp}, i.e. $\psi(x)=\sum_{i=1}^{M} w_i \phi(x,\mu_{i},\Sigma_{i})$\label{def:gauss_rbf_sum}
for some $M\in \Ne$ and $w\in\mathbb{R}^{M}$ with $\Sigma_i\succcurlyeq 0$ and $\mu_i\in\Re^n$ for $i=1,\dots,M$ where $n$ denotes the dimension of $\Xs$.
\end{defn}

In what follows, we derive conditions by which one can choose $\{w_i,\mu_i,\Sigma_i\}_{i=1}^M$ to upper bound the reach-avoid value function. We restrict ourselves to quadratic sets, defined as:
\begin{defn}[Quadratic set in $\mathbb{R}^n$] A quadratic set is any set $A \subset \mathbb{R}^n$ that can be written as 
\begin{align*}
A =\bracks*{x\in \mathbb{R}^n \Big| \begin{bmatrix}x\\1 \end{bmatrix}^\top A_j \begin{bmatrix}x\\1 \end{bmatrix}\geq 0,\ \forall  j\in\{1,\dots,N\}}
\end{align*}
for some $N\in \Ne_+$ and symmetric matrices $A_j\in\mathbb{R}^{n+1}$.
\label{def:quad_set}
\end{defn}

Sets satisfying the condition in Definition \ref{def:quad_set} include intersections of ellipsoids and half spaces that can be used to approximate convex sets arbitrarily closely. However, quadratic sets are not a subset of convex sets since one can, for example, express the difference between ellipsoids and half spaces as a quadratic set which implies that a wide range of non-convex sets satisfy Definition \ref{def:quad_set}. Given two function elements in $\F^M$, the following lemma establishes sufficient conditions on their weights, means and variances such that one uniformly dominates the other on a given quadratic set:

%To the best of our knowledge, however, there are no results to support this claim.

%%%%%%%%%%%%%%% RBF dominance %%%%%%%%%%%%%%%
\begin{lem}[RBF dominance]
\label{res:rbf_dom}
Let $\{\hat{w}_i,\hat{\mu}_i,\hat{\Sigma}_i\}_{i=1}^M$ and $\{{w}_i,{\mu}_i,{\Sigma}_i\}_{i=1}^M$ denote the weights, means and variances of two Gaussian RBF sum functions $\hat{\psi},{\psi}\in \F^M$ and assume the weights $\{w_i\}_{i=1}^M$ are non-negative. Consider a quadratic set $A\subset \mathbb{R}^n$ described by $N$ inequalities. We have that $\hat{\psi}(x)\geq \psi(x),\forall x\in A$  if for all $i=1,\dots,M$, $\exists \tau_i\geq_{\mathbb{R}^N} 0, \ \text{such that} \ X_i\succcurlyeq 0$
%\begin{align*}
%\parens*{\hat{\mu}_i,\hat{\Sigma}_i,\hat{w}_i,{\mu}_i,{\Sigma}_i,{w}_i,\tau_i}
%\exists \tau_i\geq 0, \ \text{such that} \ X_i\succcurlyeq 0
%\end{align*}
%for matrices $X_i$ defined as
\begin{align*}
X_i=
\begingroup
\renewcommand*{\arraystretch}{1.5}
\begin{bmatrix}
Q_{\hat{\Sigma}_i}^{-1}& Q_{\hat{\mu}_i} \\
Q_{\hat{\mu}_i}^\top & Q_i+Q_{\mu_i,\Sigma_i}-\sum_{j=1}^{N}\tau_{i,j} {A_j}
\end{bmatrix}
\endgroup
\end{align*}
where
\begin{align*}
Q_i= \begin{bmatrix} \mathbf{0} & \mathbf{0}\\ \mathbf{0} & 2\log\left(\tfrac{\hat{w}_i\sqrt{|\Sigma_i|}}{w_i\sqrt{|\hat{\Sigma}_i|}}\right) \end{bmatrix},\
Q_{\hat{\mu}_i}= \begin{bmatrix} \mathbf{0} & \mathbf{0} \\ \textbf{\textsc{I}} &-\hat{\mu}_i\end{bmatrix}, \  Q_{\hat{\Sigma}_i}=   \begin{bmatrix} \textbf{\textsc{I}} & \mathbf{0} \\ \mathbf{0} &{\hat{\Sigma}_i}^{-1}\end{bmatrix},\
Q_{\mu_i,\Sigma_i}= \begin{bmatrix} \Sigma^{-1}_i & -\Sigma_i^{-1}\mu_i\\-\mu_i^\top\Sigma_i^{-1} & \mu_i^\top \Sigma_i^{-1} \mu_i\end{bmatrix}
\end{align*}
and  $\mathbf{0}, \textbf{\textsc{I}}$ denote zero and identity matrices of appropriate dimensions. The terms $\tau_i\geq_{\mathbb{R}^N}0$ denote element-wise non-negativity where each $\tau_i$ is an element of space $\mathbb{R}^N$.
\end{lem}
\begin{pf}
Notice that a sufficient condition to ensure that $\hat{\psi}(x)\geq\psi(x)$, $\forall x\in A$ is that $\hat{w}_i \phi(x,\hat{\mu}_i,\hat{\Sigma}_i)\geq w_i\phi(x,{\mu}_i,{\Sigma}_i),\ \forall x\in A,\ i=1,\dots,M$.
Consider $x\in \mathbb{R}^{n}$ and let $z=\begin{bmatrix}x\\1\end{bmatrix}\in\mathbb{R}^{n+1}$. Since all $w_i$ are assumed non-negative, the weights $\hat{w}_i$ also have to be non-negative and we can take the natural logarithm on the above relation to get $\log\left(\hat{w}_i \phi(x,\hat{\mu}_i,\hat{\Sigma}_i)\right)-\log\left( w_i\phi(x,{\mu}_i,{\Sigma}_i)\right)\geq 0, \ \forall x\in A$ if and only if 
\begin{align}
\label{eq:reform_constr}
z^\top {A}_j z\geq 0,\ j=1,\dots,N \implies z^\top \left( Q_{\mu_i,\Sigma_i}-Q_{\hat{\mu}_i,\hat{\Sigma}_i}+ Q_i \right)z &\geq 0
\end{align}
%\begin{align*}
%\log\left(\hat{w}_i \phi_i(x,\hat{\mu}_i,\hat{\Sigma}_i)\right)-\log\left( w_i\phi_i(x,{\mu}_i,{\Sigma}_i)\right)&\geq 0, \ \forall x\in A \\
%\iff \\
%z^\top {R}_j z\geq 0,\ j=\{1,\dots,N\} \implies z^\top \left( Q_{\mu_i,\Sigma_i}-Q_{\hat{\mu}_i,\hat{\Sigma}_i}+ Q_i \right)z &\geq 0\\
%\end{align*}
where $Q_{\mu_i,\Sigma_i}$, $Q_{\hat{\mu}_i,\hat{\Sigma}_i}$ and $Q_i$ are defined as
\begin{align*}
Q_{\mu_i,\Sigma_i}=\bM \Sigma^{-1}_i & -\Sigma_i^{-1}\mu_i\\-\mu_i^\top\Sigma_i^{-1} & \mu_i^\top \Sigma_i^{-1} \mu_i\eM,\  
Q_{\hat{\mu}_i,\hat{\Sigma}_i}= \bM \hat{\Sigma}^{-1}_i & -\hat{\Sigma}_i^{-1}\hat{\mu}_i\\-\hat{\mu}_i^\top\hat{\Sigma}_i^{-1} & \hat{\mu}_i^\top \hat{\Sigma}_i^{-1} \hat{\mu}_i\eM,\
Q_i=\bM \mathbf{0} & \mathbf{0} \\ \mathbf{0} & 2\log\left(\tfrac{\hat{w}\sqrt{|\Sigma_i|}}{w\sqrt{\hat{\Sigma}}_i}\right) \eM 
\end{align*}
Applying the S-procedure to the condition in \eqref{eq:reform_constr}, we have that if 
\begin{align*}
\ \exists \tau_{i,j} \geq 0 \ \sthat \ Q_{\mu_i,\Sigma_i}-Q_{\hat{\mu}_i,\hat{\Sigma}_i}+ Q_i-\sum_{j=1}^N\tau_{i,j} {A}_j \succcurlyeq 0
\end{align*}
then
\begin{align*}
z^\top {A}_j z\geq 0,\ j=1,\dots,N \implies z^\top \left( Q_{\mu_i,\Sigma_i}-Q_{\hat{\mu}_i,\hat{\Sigma}_i}+ Q_i \right)z &\geq 0.
\end{align*}
Consider now the transformation $Q_{\hat{\mu}_i,\hat{\Sigma}_i}= Q_{\hat{\mu}_i}^\top Q_{\hat{\Sigma}_i}Q_{\hat{\mu}_i}$ with 
%\begin{align*}
$Q_{\hat{\mu}_i} = \bM \mathbf{0} & \mathbf{0} \\ \textbf{\textsc{I}}& -\hat{\mu}_i \eM, $ 
$Q_{\hat{\Sigma}_i}= \bM  \textbf{\textsc{I}} & \mathbf{0} \\ \mathbf{0} & \hat{\Sigma}^{-1}_i\eM $.
%\end{align*}
Notice that $Q_{\mu_i,\Sigma_i}-Q_{\hat{\mu}_i}^\top Q_{\hat{\Sigma}_i}Q_{\hat{\mu}_i}+ Q_i-\sum_{j=1}^N\tau_{i,j} A_j$ is the Schur complement of the matrix 
\begin{align*}
X_i=\bM Q_{\hat{\Sigma}_i}^{-1}  & Q_{\hat{\mu}_i}\\ Q_{\hat{\mu}_i}^\top & Q_{\mu_i,\Sigma_i}+ Q_i-\sum_{j=1}^N \tau_{i,j} A_j\eM 
\end{align*}
and since $ Q_{\hat{\Sigma}_i}^{-1}$ is positive definite by definition ($\hat{\Sigma}_i$ is the covariance matrix of a Gaussian RBF), the condition $Q_{\mu_i,\Sigma_i}-Q_{\hat{\mu}_i}^\top Q_{\hat{\Sigma}_i}Q_{\hat{\mu}_i}+ Q_i-\sum_{j=1}^N \tau_{i,j} {A}_j \succcurlyeq 0$ is equivalent to $X_i\succcurlyeq 0$.
\end{pf}
Notice that there may be cases where multipliers $\tau_i$ do not exist unless additional boundedness assumptions are imposed on the set $A$. In what follows we restrict the structure of the reach-avoid problem and illustrate how one can express the right-hand-side of the constraints in \eqref{eq:semiinfLP_chap_rob} as a sum of Gaussian RBFs. We then use the result of Lemma \ref{res:rbf_dom} to construct upper bound approximations for the reach-avoid value functions. %Note that the cost function used in Theorem \ref{res:rbf_min} is not equivalent to the one in \eqref{eq:semiinfLP_chap_rob} and can be chosen differently, affecting the quality of the bound on $\bar{\Xs}$. The constraints  

\subsection{Upper bound for the reach-avoid value functions}
\label{sec:upper_bound_analysis}
We impose the following assumptions in order to express the optimization problem in \eqref{eq:semiinfLP_chap_rob} into the framework of Gaussian RBFs. The process involves two steps: we first derive a general result concerning the expectation of a Gaussian RBF sum applied to a random variable whose density function is a Gaussian RBF sum and then apply it to the constraints in \eqref{eq:semiinfLP_chap_rob}. We make use of the  following:

%In order to enforce the upper bounding approximation to be a Gaussian RBF sum, we use upper bounds for the indicator functions $\mathds{1}_{K}$ and $\mathds{1}_{\bar{\Xs}}$ which we assume are computed a-priori.
%\nikos{Some introductory statements needed. Explain how the subsection evolves basically.}
\begin{lem}
Given a bounded set $A\subset \Re^n$, there always exists a Gaussian RBF sum $\hat{p}_A$ such that $\hat{p}_A(x)\geq \mathds{1}_A(x)$ for all $x\in\Re^n$.
\label{lem:bounded}
\end{lem}

The proof of Lemma \ref{lem:bounded} is omitted since it is trivial to show that even a single Gaussian RBF in the form of \eqref{eq:basisf_sdp} is enough to upper bound the indicator function of a bounded set by simply taking a large enough scaling weight.

%%%%%%%% Assumptions %%%%%%%%%%%%%%%
%\begin{assumption}
%\label{assum:indicator_functions}
%The target set $K$ is bounded.
%There exists a Gaussian RBF sum  $\hat{p}_K$ such that $\hat{p}_K(x)\geq \mathds{1}_K(x)$ for every $x\in \Xs$.
%\end{assumption}
%Given a bounded set $K$, there always exists a Gaussian RBF sum \hat{p}to upper bound 

% is not restrictive, provided the set $K$ is bounded and can be enclosed within a known ellipse $\m{E}$. Since the level-sets of a Gaussian RBFs are ellipses, it is straightforward to construct a Gaussian RBF that is above 1 on the whole of $\m{E}$.
%
 %The indicator function of a compact set is (in general) discontinuous but can be trivially upper bounded by a piecewise linear function which is continuous on a compact set enclosing $K$. According to the results in \cite{hartman1990layered,sandberg2001gaussian,park1991universal,cybenko1989approximation}, it is always possible to construct a Gaussian RBF sum that approximates a continuous function arbitrarily close on a given compact set.

\begin{assumption}
\label{assum:kernel} 
The stochastic kernel $Q(\cdot|x,u)$ is defined through a Gaussian RBF sum density function of the form  $\sum_{j=1}^J w_j \phi(y,\mu_j(x,u),\Sigma_j)$ for $J\in \Ne_+$ and functions $\mu_j$ affine in $(x,u)$.
%i.e. $G_(\dot,\sum_{j=1}^J \alpha_j\mathcal{N}(\mu_j,\Sigma_j)$ with known covariance matrices $\Sigma_j$, centers $\mu_j$ and weights $\alpha_j$ such that $\sum_{j=1}^J\alpha_j=1$ for a finite $J\in\mathbb{N}_+$. Note that the mixture can be different for each $k=\{0,\dots,T-1\}$.
\end{assumption}

The condition imposed on the transition kernels by Assumption \ref{assum:kernel} holds for a wide range of stochastic dynamical systems. For example, linear systems affected by general additive Gaussian mixture noise satisfy Assumption \ref{assum:kernel} and can be used as a basis to approximate the kernel of more general non-linear systems. 
%\begin{assumption}
%\label{assum:measure} 
%The state-relevance measure $\nu$ can be written as a sum of Gaussian RBFs of the form  $\sum_{j=1}^J w_j \phi(y,\nu_j,N_j)$ for some $J\in \Ne_+$.
%%i.e. $G_(\dot,\sum_{j=1}^J \alpha_j\mathcal{N}(\mu_j,\Sigma_j)$ with known covariance matrices $\Sigma_j$, centers $\mu_j$ and weights $\alpha_j$ such that $\sum_{j=1}^J\alpha_j=1$ for a finite $J\in\mathbb{N}_+$. Note that the mixture can be different for each $k=\{0,\dots,T-1\}$.
%\end{assumption}

\begin{lem}[RBF expectation]
\label{res:rbf_expect} Let $y\in\Re^n$ be a random variable distributed according to a Gaussian RBF sum density function of the form $\sum_{j=1}^J \bar{w}_j\phi(x,\bar{\mu}_j,\bar{\Sigma}_j)$ for some $J\in \mathbb{N}_+$ and $x\in\Re^n$. Consider a Gaussian RBF sum function $g\in\F^M$ defined by $g(x)=\sum_{i=1}^M w_i \phi(x,\mu_i,\Sigma_i)$ for any $x\in\Re^n$. The expected value of the random variable $g(y)$ is a Gaussian RBF sum given by $\sum_{i=1}^M \sum_{j=1}^J w_i\bar{w}_j\phi(\bar{\mu}_j,\mu_i,\Sigma_i+\bar{\Sigma}_j)$
\end{lem}
\begin{pf}
Consider the function $g(x)=\sum_{i=1}^M w_i\phi(x,\mu_i,\Sigma_i)$ and a random variable $y\in\Re^n$ with density function given by $\sum_{j=1}^J \bar{w}_j\phi(x,\bar{\mu}_j,\bar{\Sigma}_j) $, $\ x\in\Re^n$ . We have that
\begin{align*}
\expect [g(y)]&=\int_{\mathbb{R}^n} \sum_{i=1}^Mw_i\phi(x,\mu_i,\Sigma_i)\sum_{j=1}^J \bar{w}_j\phi(x,\bar{\mu}_j,\bar{\Sigma}_j) \dx\\
&=\int_{\mathbb{R}^n} \sum_{i=1}^M \sum_{j=1}^J w_i \bar{w}_j\phi(x,\mu_i,\Sigma_i)\phi(x,\bar{\mu}_j,\bar{\Sigma}_j) \dx\\
&=\sum_{i=1}^M \sum_{j=1}^J w_i\bar{w}_j\int_{\mathbb{R}^n} \phi(x,\mu_i,\Sigma_i)\phi(x,\bar{\mu}_j,\bar{\Sigma}_j) \dx\\
&=\sum_{i=1}^M \sum_{j=1}^J w_i\bar{w}_j\phi(\bar{\mu}_j,\mu_i,\Sigma_i+\bar{\Sigma}_j)
\end{align*}
 The last equality follows from the fact that the product of two functions in the form of (\ref{eq:basisf_sdp}) is proportional to another function of the same form \cite[Section 2]{hartman1990layered}. With
\begin{align*}
\phi(x,\mu_i,\Sigma_i)\phi(x,\bar{\mu}_j,\bar{\Sigma}_j)=\phi(\mu_i,\bar{\mu}_j,\Sigma_i+\bar{\Sigma}_j)\phi(x,\bar{\mu},\bar{\Sigma}),
\end{align*}
$\bar{\mu}=\left(\Sigma_i^{-1}+\bar{\Sigma}_j^{-1}\right)^{-1}\left(\Sigma_i^{-1}\mu_i+\bar{\Sigma}_j^{-1}\bar{\mu}_j\right)$, and $\bar{\Sigma}=\left(\Sigma_i^{-1}+\bar{\Sigma}_j^{-1}\right)^{-1}$. The result follows as the proportionality constant $\phi(\mu_i,\bar{\mu}_j,\Sigma_i+\bar{\Sigma}_j)$ is independent of $x$ and the function $\phi(x,\bar{\mu},\bar{\Sigma})$ integrates to 1 over $\mathbb{R}^n$. 
\end{pf}

%Define the approximation operator $\hat{\mathcal{T}}[V](x) := \hat{p}_K(x)+\hat{p}_{\bar{\Xs}}(x)\int_{\mathcal{X}}{V}(y)G_k(dy|x)$ where $\hat{p}_K,\hat{p}_{\bar{\Xs}}$ are defined by Assumption \ref{assum:indicator functions}. It clearly  holds that $\hat{\mathcal{T}}[V](x)\geq \mathcal{T}[V](x)$ for all $x\in\R^n$ and $V\in \F$. We will use $\hat{\m{T}}$ to compute an upper bound to the solution function of \eqref{eq:inflp_ub}. Note that using the same machinery one can compute a lower bound to the solution of \eqref{eq:inflp_lb} by imposing Assumption \ref{assum:indicator functions} with reversed inequalities.

\begin{prop}[Upper bound] 
Consider Lemma \ref{lem:bounded} and Assumption \ref{assum:kernel}. Let $\hat{V}^*_T(x)=\hat{p}_K(x)$ where $\hat{p}_K$ upper bounds $\mathds{1}_K$. The following recursion of optimization problems provides an upper bound on the reach-avoid value function for each $k=T-1,\dots,0$, for any choice of non-negative measure $\nu$.
\begin{align}
\label{opt:semiinf_gauss}
\begin{split}
\hat{V}^*_{k}=\arg\min_{\hat{V}_k\in\tilde{\F}}\quad& \int_{\Xs}\hat{V}_k(x)\nu(\dx)\\
\st \quad &\hat{V}_k(x)\geq\int_{\Xs}\hat{V}_{k+1}^*(z)Q(\dz|x,u), \quad \forall (x,u)\in \bar{\Xs}\times \Us\\
\quad &\hat{V}_k(x)\geq 1\quad \forall x \in K.
\end{split}
\end{align}
\end{prop}

\begin{pf}
We use induction. The reach-avoid value function at time $T-1$ is defined by 
\begin{align*}
V^*_{T-1}(x)=\max_{u\in \Us}\left\{\mathds{1}_K(x)+\mathds{1}_{\bar{\Xs}}(x)\int_{\Xs}{V}^*_T(z)Q(\dz|x,u)\right\}, \forall x\in\Xs.
\end{align*} 
By definition we have that $\hat{V}_T^*(x)=\hat{p}_K(x)\geq \mathds{1}_K(x)=:V^*_T(x)$ for all $x\in \Xs$. We also have that for all $(x,u)\in\Xs\times\Us$, $\int_{\Xs}\hat{V}^*_T(z)Q(\dz|x,u)\geq \int_{\Xs}V^*_T(z)Q(\dz|x,u)$. Hence, for any feasible point $\hat{V}_{T-1}$ of \eqref{opt:semiinf_gauss}, we have that $\hat{V}_{T-1}(x)\geq \mathds{1}_{\bar{\Xs}}(x)\int_{\Xs}\hat{V}^*_T(z)Q(\dz|x,u)\geq \mathds{1}_{\bar{\Xs}}(x)\int_{\Xs}{V}^*_T(z)Q(\dz|x,u)$ for all $(x,u)\in\Xs\times\Us$. Moreover, the constraint $\hat{V}_{T-1}(x)\geq 1,\ \forall x \in K$ ensures that all feasible functions satisfy $\hat{V}_{T-1}(x)\geq \mathds{1}_K(x)$ for all $x\in\Xs$. We conclude that for any feasible function in \eqref{opt:semiinf_gauss} it holds that $\hat{V}_{T-1}(x)\geq \mathds{1}_K(x)+\mathds{1}_{\bar{\Xs}}(x)\int_{\Xs}{V}^*_T(z)Q(\dz|x,u), \ \forall (x,u)\in \bar{\Xs}\times \Us$ and it thus upper bounds $V^*_{T-1}(x)$. Using induction the same result follows for all $k=T-2,\dots,0$.
\end{pf}
The final step of the process is to derive a tractable convex optimization program to construct a feasible solution to the problem in \eqref{opt:semiinf_gauss}. We can then use it recursively to construct upper bounds to the value function for all $k=0,\dots,T-1$. First notice that from Lemma \ref{lem:bounded} and Assumption and \ref{assum:kernel} both $\hat{V}_{T}^*$ and $Q(\cdot|x,u)$ are Gaussian RBF sums and the right-hand-side of the first constraint in \eqref{opt:semiinf_gauss} is an RBF expectation that according to Lemma \ref{res:rbf_expect} can be written as another Gaussian RBF sum. As a result, by choosing the number of basis to construct $\tilde{\F}$ according to the number of elements in the RBF expectation in the constraints of \eqref{opt:semiinf_gauss}, we can use the RBF dominance result in Lemma \ref{res:rbf_dom} to reformulate the constraints into sufficient semidefinite constraints. Notice however that depending on the number of RBF elements used to describe $Q$ and the number of RBF elements used to approximate the value function at step $k=T$, the number of element required to build $\tilde{\F}$ at every step of the recursion may increase, affecting the computational complexity of the process. The following proposition provides a tractable semidefinite program (SDP) to construct an upper bound to the reach-avoid value function.
\begin{prop}
\label{res:vf_upper}
%%%%%%%%%%%%%
Consider Lemma \ref{lem:bounded} and Assumption \ref{assum:kernel} and further assume that the sets $K\subset K^\prime$ and $S:=\bar{\Xs}\times \Us=(K^\prime\setminus K) \times \Us $ are quadratic sets in $\mathbb{R}^n$ and $\mathbb{R}^{n+m}$ respectively. Let $\hat{V}^*_{k+1}(x)=\sum_{i=1}^M w_i \phi(x,\mu_i,\Sigma_i)$ be a known Gaussian RBF sum for some $M\in \Ne_+$ for which it holds that $\hat{V}_{k+1}(x)\geq V^*_{k+1}(x)$ for all $x\in \Xs$. Consider a density function for $Q$ given by a single Gaussian RBF defined as $q(y,x,u)= \phi(y,\mu_0(x,u),\Sigma_0)$ with the function $\mu_0:\Xs\times \Us\rightarrow \Xs$ defined as $\mu_0(x,u)=Ax +Bu$ where $A,B$ are known matrices of appropriate dimensions. Let $\nu$ be a probability measure over $\Xs$ and consider the following semidefinite program:
\begin{align}
\label{opt:boundSDPrbf}
\begin{split}
%\min_{\mu_i,\Sigma_i,w_i} &\int_A \sum_{i=1}^M w_i\phi(x,\mu_i,\Sigma_i) \dx\\
%\min_{\hat{\mu}_i,\hat{\Sigma}_i,y_i,\tau_i,\rho_i}\quad &\sum_{i=1}^{M}\sum_{j=1}^J c_{i,j}\left(y_i + \tfrac{1}{2}\log|\hat{\Sigma}_i|\right)\\
%\min_{\hat{\mu}_i,\hat{\Sigma}_i,y_i,\tau_i,\rho_i,\gamma_i}\quad &\sum_{i=1}^{M} y_i + \trace(\hat{\Sigma})_i+N) -\gamma_i\\
\min_{\hat{\mu}_i,\hat{\Sigma}_i,y_i,\tau_i,\rho_i}\quad &\sum_{i=1}^{M} y_i + \tfrac{1}{2}\trace(\hat{\Sigma}_i) \\
%\min_{\hat{\mu}_i,\hat{\Sigma}_i,y_i,\tau_i,\rho_i,\gamma_i}\quad &\sum_{i=1}^{M}y_i + \trace \\
\st \quad &\begin{rcases}
X_i\succcurlyeq 0\\
\bar{X}_i\succcurlyeq 0\\
%C_i \succcurlyeq 0\\
\tau_i,\rho_i\geq_{\mathbb{R}^N} 0 \\
\hat{\mu}_i\in \mathbb{R}^n\\
\hat{\Sigma}_i\in \mathbf{S}_+^n
	\end{rcases}\ i=1,\dots,M
	\end{split}	
\end{align}
where 
\begin{align*}
X_i=
\begingroup
\renewcommand*{\arraystretch}{1.5}
\begin{bmatrix}
Q_{\hat{\Sigma}_i}^{-1} & Q_{\hat{\mu}_i} \\
Q_{\hat{\mu}_i}^\top & Q_i+Q_{\mu_i,\bar{\Sigma}_i}-\sum_{j=1}^{N}\tau_{i,j} {S_j}
\end{bmatrix}
\endgroup,\quad
\bar{X}_i=
\begingroup
\renewcommand*{\arraystretch}{1.5}
\begin{bmatrix}
 Q_{\hat{\Sigma}_i}^{-1} & Q_{\hat{\mu}_i} \\
Q_{\hat{\mu}_i}^\top &\bar{Q}_i-\sum_{j=1}^{N}\rho_{i,j} {K_j} 
\end{bmatrix}
\endgroup
%C_i=\begingroup
%\renewcommand*{\arraystretch}{1.5}
%\begin{bmatrix}
%\gamma_i  & \nu-\hat{\mu}_i \\ (\nu-\hat{\mu}_i)^\top & \hat{\Sigma}_i+N
%\end{bmatrix}
%\endgroup
\end{align*}
with
\begin{align*}
%Q_i&= \begin{bmatrix} \mathbf{0} &\mathbf{0} &\mathbf{0}\\\mathbf{0} &\mathbf{0} &\mathbf{0}\\ \mathbf{0} & \mathbf{0}& 2\log\left(\tfrac{\hat{w}_i\sqrt{|\Sigma_i|}}{w_i\sqrt{|\hat{\Sigma}_i|}}\right) \end{bmatrix},\
Q_i&= \begin{bmatrix} \mathbf{0} &\mathbf{0} &\mathbf{0}\\\mathbf{0} &\mathbf{0} &\mathbf{0}\\ \mathbf{0} & \mathbf{0}& 2\left(y_i+\log\left(\tfrac{\sqrt{|\bar{\Sigma}_i|}}{w_i}\right)\right)\end{bmatrix},\
Q_{\hat{\mu}_i}= \begin{bmatrix} \mathbf{0} & \mathbf{0} & \mathbf{0}\\\mathbf{0} &\mathbf{0} &\mathbf{0} \\ \textbf{\textsc{I}}&\mathbf{0} &-\hat{\mu}_i\end{bmatrix}, \  
%%%%%%%%%%%%%%%%%%%%%%%
Q_{\hat{\Sigma}_i}=   \begin{bmatrix} \textbf{\textsc{I}} & \mathbf{0}&\mathbf{0}\\\mathbf{0} &\textbf{\textsc{I}} &\mathbf{0} \\ \mathbf{0} &\mathbf{0}&{\hat{\Sigma}_i}^{-1}\end{bmatrix},\\
%%%%%%%%%%%%%%%%%%%%%%%
Q_{\mu_i,\bar{\Sigma}_i}&= \begin{bmatrix} A^\top\bar{\Sigma}_iA & A^\top \bar{\Sigma}_iB &-A^\top\bar{\Sigma}_i\mu_i\\B^\top\bar{\Sigma}_iA &B^\top\bar{\Sigma}_iB & -B^\top \bar{\Sigma}_i\mu_i\\-\mu_i^\top\bar{\Sigma}_i A &-\mu_i^\top \bar{\Sigma}_i B& \mu_i^\top \bar{\Sigma}_i\mu_i\end{bmatrix},\
%%%%%%%%%%%%%%%%%%%%%%%
%\bar{Q}_i= \begin{bmatrix} \mathbf{0} & \mathbf{0}\\ \mathbf{0} &2\log\left(\tfrac{(M+J)\hat{w}_i}{\sqrt{(2\pi)^n|\hat{\Sigma}_i|}}\right) \end{bmatrix},\
\bar{Q}_i= \begin{bmatrix} \mathbf{0} &\mathbf{0} & \mathbf{0}\\\mathbf{0} &\mathbf{0} & \mathbf{0}\\ \mathbf{0} &\mathbf{0}&2\left(y_i+\log\left(\tfrac{M}{\sqrt{(2\pi)^n}}\right)\right) \end{bmatrix}
%%%%%%%%%%%%%%%%%
\end{align*}
and $y_i:=\log\left(\tfrac{\hat{w}_i}{\sqrt{|\hat{\Sigma}_i|}}\right)$, $\bar{\Sigma}_i=\left({\Sigma_i}+\Sigma_0\right)^{-1}$. Let $\{\hat{\mu}^*_i,\hat{\Sigma}^*_i,y^*_i,\tau^*_i,\rho^*_i\}_{i=1}^M$ denote the optimal solution of \eqref{opt:boundSDPrbf} and set $\hat{w}_i^*=e^{y_i^*}\sqrt{|\hat{\Sigma}^*_i|}$. The function constructed by $\hat{V}^*(x)=\sum_{i=1}^M \hat{w}^*_i\phi\left(x,\hat{\mu}_i^*,\hat{\Sigma}^*_i\right)$ is a feasible solution for \eqref{opt:semiinf_gauss}.
\end{prop}
\begin{pf}
First we show how the constraints of problem \eqref{opt:boundSDPrbf} are sufficient for the constraints of \eqref{opt:semiinf_gauss}. Given the format of functions $\hat{V}^*_{k+1}$ and $q$ we have from Lemma \ref{res:rbf_expect} that
\begin{align*}
\int_{\Xs}\hat{V}^*_T(z)Q(\dz|x,u) &=\int_{\Xs} \sum_{i=1}^M w_i \phi(z,\mu_i,\Sigma_i)q(z,x,u)\dz\\
&=\int_{\Xs} \sum_{i=1}^M w_i \phi(z,\mu_i,\Sigma_i)\phi(z,\mu_0(x,u),\Sigma_0)\dz\\
&= \sum_{i=1}^M w_i\phi\left(\mu_0(x,u),\mu_i,\Sigma_i+\Sigma_0\right).
\end{align*}
Since $\hat{V}(x)=\sum_{i=1}^M \hat{w_i}\phi(x,\hat{\mu}_i,\hat{\Sigma}_i)$, the constraint $\hat{V}(x)\geq\int_{\Xs}\hat{V}^*_T(z)Q(\dz|x,u)$ in \eqref{opt:semiinf_gauss} is equivalent to the RBF dominance constraint 
\begin{align*}
\sum_{i=1}^M \hat{w_i}\phi(x,\hat{\mu}_i,\hat{\Sigma}_i)\geq \sum_{i=1}^M w_i\phi\left(\mu_0(x,u),\mu_i,\Sigma_i+\Sigma_0\right), \quad \forall (x,u)\in S
\end{align*}
which according to Lemma \ref{res:rbf_dom} is implied by the constraints $X_i\succcurlyeq 0$ and $\tau_i \geq 0$ for $i=1,\dots,M$. The form of the matrices $Q_i, Q_{\mu_i,\Sigma_i}, Q_{\hat{\mu}_i}, Q_{\hat{\Sigma}_i}$ is equivalent to the one derived in the proof of Lemma \ref{res:rbf_dom}, considering the variable $z=\bM x\\u\\1 \eM$, the function $\mu_0(x,u)=Ax+Bu$ and the transformation $y_i=\log\left(\tfrac{\hat{w}_i}{\sqrt{|\hat{\Sigma}_i|}}\right)$. We replace the second constraint $\hat{V}(x)\geq 1,\ \forall x \in K$ by the sufficient condition 
%\begin{align*}
$\sum_{i=1}^M \hat{w_i}\phi(x,\hat{\mu}_i,\hat{\Sigma}_i)\geq \sum_{i=1}^M \tfrac{1}{M},\ \forall x\in K$
%\end{align*}
which is again an RBF dominance constraint that according to Lemma \ref{res:rbf_dom} is implied by the constraints $\bar{X}_i\succcurlyeq 0$ and $\tau_i \geq 0$ for $i=1,\dots,M$. Finally, the cost function in \eqref{opt:semiinf_gauss} reads
\begin{align*}
\int_{\Xs}\hat{V}(x)\nu(\dx)&=\int_{\Xs}\sum_{i=1}^M \hat{w_i}\phi(x,\hat{\mu}_i,\hat{\Sigma}_i)\dx=\sum_{i=1}^M \hat{w_i}=\sum_{i=1}^M e^{y_i}\sqrt{|\hat{\Sigma}_i|}
\end{align*}
which we replace with the sum of the logarithms, i.e. $\sum_{i=1}^M \log\left(e^{y_i}\sqrt{|\hat{\Sigma}_i|}\right)=\sum_{i=1}^M y_i +\tfrac{1}{2}\log(|\hat{\Sigma}_i|)$ motivated by the fact that the logarithm is a monotone function. Still, the elements $\tfrac{1}{2}\log(|\hat{\Sigma}_i|)$ are concave in $\hat{\Sigma}_i$ and the minimization is a non-convex problem. Using the fact that $\tfrac{1}{2}\log(|\hat{\Sigma}_i|)=\tfrac{1}{2}\log(\prod_{d=1}^n \lambda_d^i)=\tfrac{1}{2}\sum_{d=1}^n \log(\lambda_d^i)$ and $\trace(\hat{\Sigma}_i)=\sum_{d=1}^n \lambda_d^i$ we employ $\tfrac{1}{2}\trace(\hat{\Sigma}_i)$ as a heuristic to approximate the terms $\tfrac{1}{2}\log(|\hat{\Sigma}_i|)$. 
\end{pf}
Feasibility of the problem in Proposition \ref{res:vf_upper} is not easy to prove without additional assumptions. One way to guarantee a feasible solution is to add the constant function $c(x)=1, \forall x\in \Xs$ in the basis set used to construct $\tilde{\F}$. Notice that we have expressed the kernel $Q$ as a single Gaussian RBF to simplify the notation in \eqref{opt:boundSDPrbf}; extensions to Gaussian RBF sums with varying functions $\mu_0(x,u)$ are straightforward.

The SDP-based procedure presented in Proposition \ref{res:vf_upper} can be used recursively, starting at $k=T-1$ and $\hat{V}^*_T(x)=\hat{p}_K(x)$ for all $x\in \Xs$, to construct upper bound approximations to the value functions for every $k=0,\dots,T-1$. Note that the process introduces conservatism due to the RBF dominance constraints used in \eqref{opt:boundSDPrbf} which are only sufficient conditions for the original robust constraints in \eqref{opt:semiinf_gauss}. Moreover, the cost function in \eqref{opt:semiinf_gauss} is non-convex jointly in the weights, means and variances and the approximation used in Proposition \ref{res:vf_upper} constitutes a design choice that can affect how close the bounding functions are to the optimal ones on different areas of the state space. Finally, the state-relevance measure is assumed to be a probability measure on $\Xs$ but other choices allowing calculation of the integral $\int_{\Xs}\hat{V}(x)\nu(\dx)$ as a function of the RBF parameters can be used (for example measures supported on hyper-rectangular sets - see \cite{kariotoglou2014adp}). Algorithm \ref{algo:rbf_robust} summarizes the process of using Proposition \ref{res:vf_upper} recursively to compute bounds down to time $k=0$. Each SDP in the process has a total of $3M$ semidefinite constraints where the ones corresponding to $X_i$ are of dimension $n+m+1$, the ones corresponding to $\bar{X}_i$ are of dimension $n+1$ and the ones corresponding to $\hat{\Sigma}_i$ are of dimension $n$. If more than one Gaussian RBF elements are used to describe the process kernel $Q$, the total number of semidefinite constraints grows combinatorially in the elements of $Q$ and $M$ and linearly in the horizon length (see Lemma \ref{res:rbf_expect}).
\begin{algorithm}[t]
\caption{Approximate value function}
\label{algo:rbf_robust}
\begin{algorithmic}[1]
\STATE \textbf{Input Data:} 
\begin{itemize}
\small
  \setlength{\itemsep}{0.5pt}
	\item Quadratic set constraints for $\bar{\Xs}$, $K$ and $\Us$.
	\item Reach-avoid time horizon $T$.
	\item Center and variance of the MDP kernel $Q$.
		\item Gaussian RBF sum $\hat{p}_K$ for which it holds that $\hat{p}_K(x)\geq \mathds{1}_K(x)$ for all $x\in \Xs$.
\end{itemize}
\STATE \textbf{Design parameters:} 
\begin{itemize}
\small
  \setlength{\itemsep}{0.5pt}
	\item State-relevance measure $\nu$.
	\item Approximation of the cost function $\int_{\Xs}\hat{V}(x)\nu(\dx)$.		
	\end{itemize}
%\STATE Choose the required violation levels $\bracks*{\varepsilon_0,\dots,\varepsilon_{T-1}}$ and confidences $\bracks*{\beta_0,\dots,\beta_{T-1}}$. 
\STATE Initialize $\hat{V}^*_T(x)\leftarrow \hat{p}_{K}(x)$.
\FOR {$k=T-1$ \TO $k=0$ }
	\STATE Extract the weights, means and variances $\{w_i,\mu_i,\Sigma_i\}_{i=1}^M$ from the known function $\hat{V}_{k+1}^*$.
	\STATE Solve the optimization program in \eqref{opt:boundSDPrbf} to obtain $\{\hat{\mu}^*_i,\hat{\Sigma}^*_i,y^*_i\}_{i=1}^M$.
	\STATE Compute $\hat{w}^*_i=e^{y_i^*}\sqrt{|\hat{\Sigma}^*_i|}$ for $i=1,\dots, M$.
	\STATE Set $\hat{V}^*_k(x)=\sum_{i=1}^M\hat{w}_i^*\phi(x,\hat{\mu}_i^*,\hat{\Sigma}^*_i)$.
\ENDFOR
\end{algorithmic}
\end{algorithm}

\subsection{Approximate control policy}
Given the value function $V^*_k$ at each $k=\{0,\dots,T-1\}$, the problem in \eqref{eq:optpol} implicitly defines the optimal reach-avoid control policy at time $k$. Using the sequence of approximate value functions $\hat{V}^*_0,\dots,\hat{V}^*_{T-1}$constructed using Proposition \ref{res:vf_upper} recursively, one can compute approximate control policies $\hat{\pi}_k$ at each step $k=\{0,\dots,T-1\}$. If the kernel $Q$ has a Gaussian RBF sum density function $q(y,x,u)=\phi(y,Ax+Bu,\Sigma)$, as in Proposition \ref{res:vf_upper}, given $\hat{V}^*_{k+1}(x)=\sum_{i=1}^M w_i \phi(x,\mu_i,\Sigma_i)$, we define the SDP-based approximate control policy at each time step $k=\{0,\dots,T-1\}$ as:
\begin{align}
\begin{split}
\hat{\pi}_k(x)&=\arg\max_{u\in \mathcal{U}}\left\{\mathds{1}_K(x)+\mathds{1}_{\bar{\Xs}}(x)\int_{\Xs}{\hat{V}^*_{k+1}(z)Q(\dz|x,u)}\right\}\\
&=\arg\max_{u\in \mathcal{U}}\int_{\Xs}{\hat{V}_{k+1}^*(z)Q(\dz|x,u)}\\
&\overset{(\star)}=\arg\max_{u\in \mathcal{U}}\sum_{i=1}^M w_i\phi (Ax+Bu,\mu_i,\Sigma_i+\Sigma)
\end{split}
\label{eq:approx_pol_sdp}
\end{align}
where we have used Lemma \ref{res:rbf_expect} to compute the expectation in $(\star)$. Despite the fact that the optimization problem in \eqref{eq:approx_pol_sdp} is non-convex, standard gradient based algorithms can be employed to obtain a local solution. In particular, the cost function is by construction smooth with respect to $u$ for a fixed $x\in\bar{\Xs}$ and the gradient and Hessian functions can be analytically computed. Moreover, the decision space $\Us$ is typically low dimensional (in most mechanical systems for example $\dim{\Us}<\dim{\Xs})$ and mature software is available \cite{waechter2009introduction} to compute locally optimal solutions. The process of calculating a control input at time $k$ for a fixed state $x_k$ is summarized in Algorithm \ref{algo:controlalgo_sdp}. Alternative approaches include using randomized techniques similar to the approaches in \cite{calafiore2006scenario,campi2008exact,petretti2014approximatelinear,deori2014buildingrandomized}, general purpose non-linear programming solvers (e.g. interior point method) or directly gridding the control space, if computationally feasible, as done in the second example of the next section. 
\begin{algorithm}[t]
\caption{Approximate controller}
\label{algo:controlalgo_sdp}
\begin{algorithmic}[1]
\FOR {$k\in \bracks*{0,\dots,T-1}$}
\STATE Measure the system state $x_k$.
\STATE Compute the gradient and hessian of $\sum_{i=1}^M w_i\phi (Ax_k+Bu,\mu_i,\Sigma_i+\Sigma)$ with respect to $u$.
\STATE Solve the optimization problem in \eqref{eq:approx_pol_sdp} to obtain the control input $\hat{\pi}_k(x_k)$.
\STATE Apply the calculated control input to the system.
\ENDFOR
\end{algorithmic}
\end{algorithm}

%In the following section, we evaluate the performance of the approximate control policies defined by $\hat{\mu}_k$ for $k=\{0,\dots,T-1\}$. Starting at some $x_0\in \bar{\Xs}$, we implement the policy and gather statistics via Monte Carlo-type stochastic simulations where we record how many times a given MDP satisfies the desired reach-avoid objective. 
\section{Numerical examples}
\label{sec:examples}
In this section we evaluate the SDP-based upper bounds calculated using the process of Algorithm \ref{algo:rbf_robust}. We first investigate the suboptimality of the bound and the performance of the associated control policies by comparing them to value functions and policies computed via space gridding. Since space gridding is limited to low dimensional state-input spaces ($\dim(\Xs\times \Us)$ around 4-5), we also evaluate the bounding value functions on a simple benchmark problem where we assume an LQG controller is a suitable heuristic to maximimize the reach-avoid probability. All simulations were carried out on an Intel Core i7 Q820 CPU clocked @ 1.73 GHz with 16GB of RAM memory, using MOSEK's SDP solver in its default settings.

\subsection{Reach-avoid probability bounds compared to space gridding}
\label{sec:compare_grid_sdp}
\begin{figure}[t]
\begin{minipage}[t]{0.45\linewidth}
\centering
\vspace{0pt}
\includegraphics[width=7cm, height=5cm]{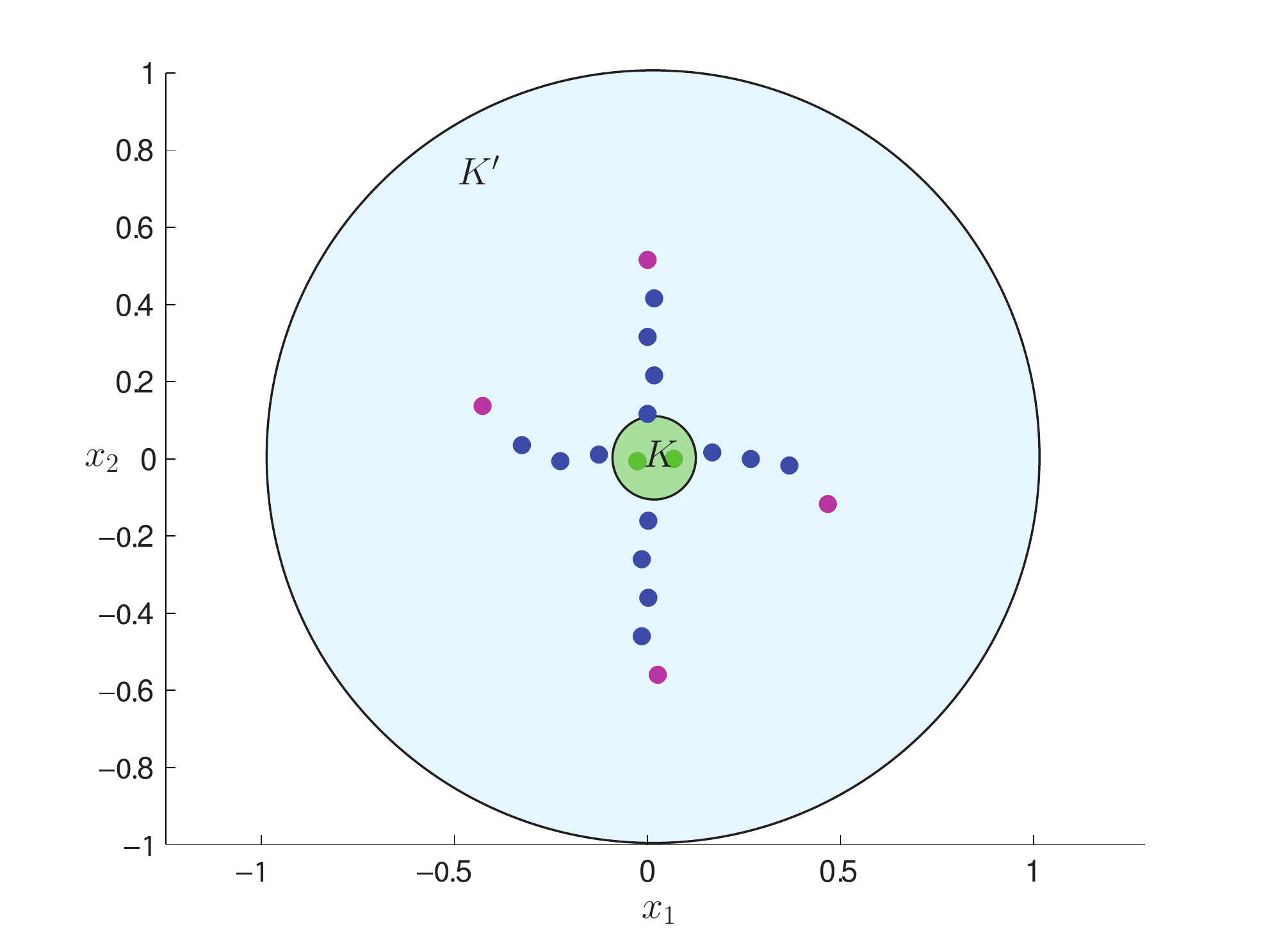}%
\caption{Reach-avoid safe and target sets with example trajectories}%
\label{fig:4Dspaces}%
\end{minipage}\hfill
\begin{minipage}[t]{0.45\linewidth}
\centering
\vspace{0pt}
\includegraphics[width=7cm, height=5cm]{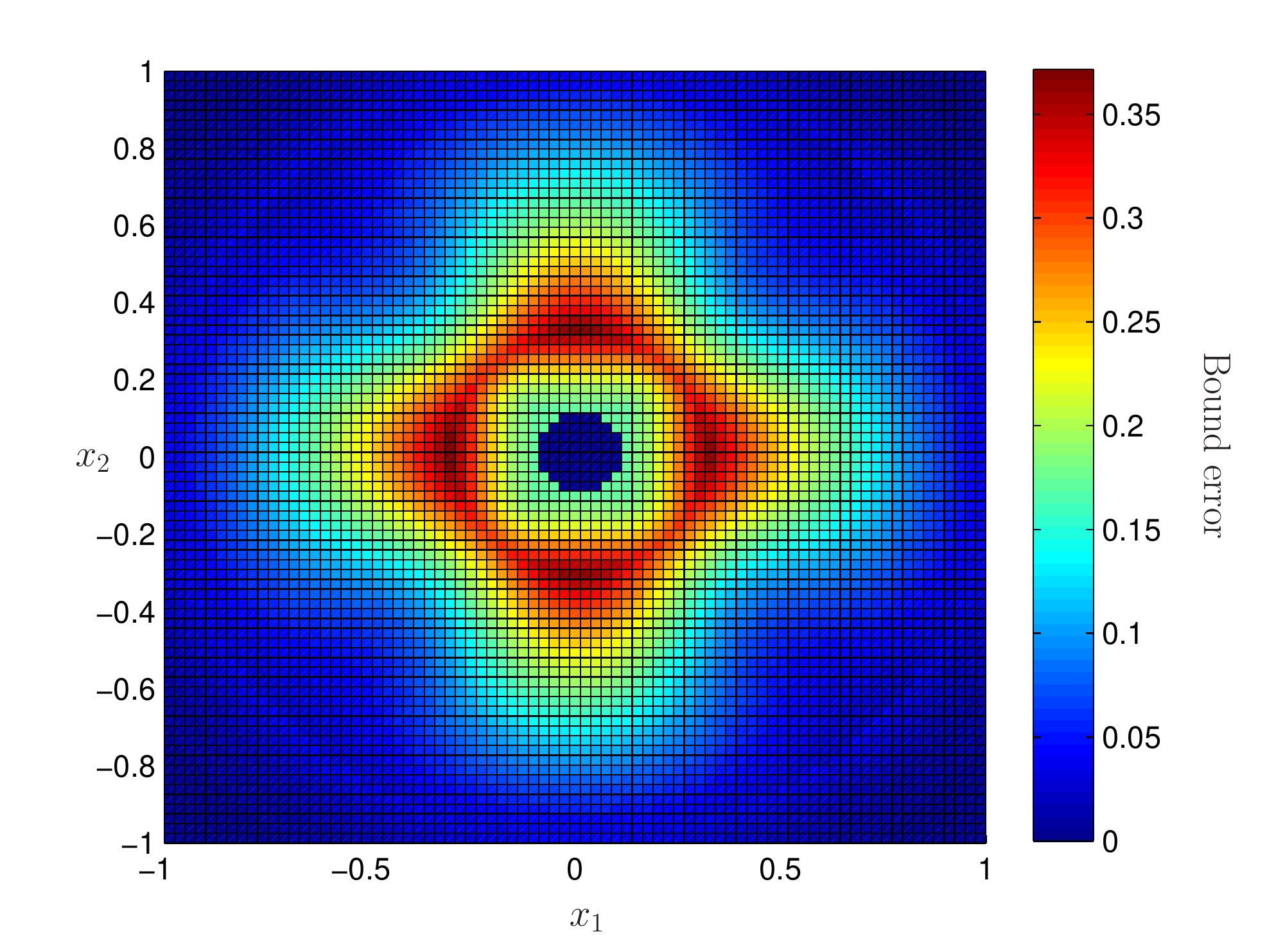}%
\caption{Value function bound error at time $k=0$}%
\label{fig:vf_sdp}%
\end{minipage}
\end{figure}

Consider the reach-avoid problem of maximizing the probability that the state of a controlled linear system subject to additive Gaussian noise reaches a target set around the origin within $T=5$ discrete time steps while staying in a safe set for all previous times. We consider systems described by the equation 
\begin{align}
	\label{eq:linsys_sdp}
		\begin{split}
x_{k+1}= Ax_k + Bu_k +\omega_k
\end{split}
\end{align} 
where for each $k=\{0,\dots,4\}$, $x_k \in \Xs= \mathbb{R}^n$, $u_k\in \Us= \{u\in \mathbb{R}^m| u^\top Q_u u\leq \rho_u^2\}$ with $\rho_u=0.1$ and each $\omega_k$ is distributed as a Gaussian random variable $\omega_k\sim\mathcal{N}\parens*{\mathbf{0}_{n\times 1},\Sigma_\omega}$ with $\Sigma_\omega \in \mb{S}_+^n$. We consider a target set $K=\{x\in \mathbb{R}^n | x^\top Q_t x\leq \rho_t^2\}$ with $\rho_t=0.1$ around the origin and a safe set $K=\{x\in \mathbb{R}^n | x^\top Q_s x\leq \rho_s^2\}$ with $\rho_s=1$ (see Figure \ref{fig:4Dspaces}). 
In order to easily change the dimensions of the problem (variables $n,m$) and keep the system behaviour between problems similar, we fix $A=\mb{1}_{n\times n}$, $B=\mb{1}_{n\times m}$ and $\Sigma_\omega=0.001\cdot \mb{1}_{n\times n}$. 

For problem dimensions up to $n=3$ and $m=2$, we compute the gridding-based value functions by discretizing the spaces $[-\rho_s,\rho_s]^n$ and $[-\rho_u,\rho_u]^m$ and solving the recursion in \eqref{eq:DP} on the resulting grid, initialized with $V^*_{T}(x)=\mathds{1}_K(x)$. Table \ref{tab:grid_comp} reports the space discretization for each problem dimension ($N_s$ corresponds to number of elements in each dimension of $[-\rho_s,\rho_s]$ and $N_u$ corresponds to number of elements in each dimension of $[-\rho_u,\rho_u]$) along with the time it took to compute a single iteration (one horizon step) of the recursion. The upper bound value functions are constructed using the steps of Algorithm \ref{algo:rbf_robust}. The transition kernel of \eqref{eq:linsys_sdp} has a Gaussian RBF sum density $x_{k+1}\sim \mathcal{N}(Ax_k+Bu_k,\Sigma_\omega)$ with center $Ax_k+Bu_k$ and covariance matrix $\Sigma_\omega$ while the sets $K$, $K^\prime$ and $\Us$ are by definition quadratic sets, satisfying the input requirements of the algorithm. As in Proposition \ref{res:vf_upper}, we choose a uniform  state-relevance measure $\nu$ supported on ${\Xs}$. The only remaining element to execute the algorithm is to construct an upper bound function $\hat{p}_K$ for the indicator function $\mathds{1}_K$. We achieve this by randomly placing Gaussian RBF elements on $K$ with diagonal covariance matrices $\Sigma_b=0.0005$, choosing their weights by means of a linear program, making sure that they upper bound $\mathds{1}_K$ on every point of the grid placed on $\Xs$ (see \cite{kariotoglou2014adp} for details). Table \ref{tab:grid_comp} reports the number of elements used to upper bound $\mathds{1}_K$ for the different cases of $\dim(\Xs)=1,2,3$.

To estimate the performance of the approximate control policy, we sample 100 initial states $x_0$ uniformly from $\bar{\Xs}$ and for each one generate 100 different noise trajectory realizations, using the grid-based and SDP-based policies computed by \eqref{eq:optpol} and Algorithm \ref{algo:controlalgo_sdp} respectively. We then count the number of trajectories that successfully complete the reach-avoid objective, i.e. reach $K$ without leaving $K^\prime$, making sure that $u\in\Us$. In Table \ref{tab:grid_comp} we denote by $\hat{\pi}-\pi^*$ the mean difference between the empirical success probability of the SDP-based controller ($\hat{\pi}$) and the empirical success probability of the controller based on the value function computed via space gridding ($\pi^*$). The value of $\hat{V}^*_0-V^*_0$ corresponds to the mean difference between the values of the upper bound ($\hat{V}^*_0$) saturated to 1, and the gridding-based value function ($V^*_0$) at time $k=0$, over the grid elements within the set $\bar{\Xs}$. Analyzing the results, we notice that the approximate control policy based on the value function upper bound yields near optimal performance even though the average absolute error is non-zero. Overall, the SDP-based controller outperforms the grid-based controller, probably due to coarse discretization of $\Xs$ and $\Us$. In terms of computation time, the SDP-based approach  significantly outperforms space gridding. The time values reported in Table \ref{tab:grid_comp} correspond to the computation time required to solve one step of the recursions; it is clear that the curse of dimensionality sets a fundamental limit to potential applications of space gridding.

\begin{table}[t]
\centering
\begin{tabular}{lllll|llll}
$n$&$m$&$N_s$&$N_u$&$M$&Grid (sec.)& SDP bound (sec.)&$\hat{V}_0^*-V^*_0$&$\hat{\pi}-\pi^*$\\
\hline
1&1&80&25&10&0.4262&0.62&0.34&0.035\\
2& 1&80&25&50&115&1.88&0.174&-0.0186\\
2& 2&80&25&50&1237 &1.83&0.472&-0.087\\
3&1 &30&15&90& 477&4.8&0.053&0.0083\\
3&2 &30&15&90&2926&4.15&0.186&0.031\\
\hline
\end{tabular}
\caption{Comparison of grid-based and upper bound value functions}
\label{tab:grid_comp}
\end{table}

\subsection{High-dimensional problems}
In this section we demonstrate the scalability of the SDP-based method by comparing the performance of the control policy computed using the value function bounds constructed with Algorithm \ref{algo:rbf_robust} to an LQG controller. Instead of using Algorithm \ref{algo:controlalgo_sdp} to generate the control policy, we discretize the control space $\Us$ with 20 points in each dimension, resulting in a $20^m$ grid over which we compute the approximate control policy. We refrain from using non-linear programming solvers in this comparison as computation times grow significantly with the dimension of $\Us$ and are impractical for Monte Carlo-type verification simulations. 

We consider the same regulation problem as in Section \ref{sec:compare_grid_sdp} but randomly choose the matrices $A\in \mathbb{R}^{n\times n}$ and $B\in\mathbb{R}^{n\times m}$ in every experiment such that the resulting systems are open loop stable and controllable under zero noise. If the system is unstable, the reach-avoid probability over a short horizon can be close to zero, especially when $\dim(\Xs) >> \dim(\Us)$. For the series of experiments carried out in this part, we upper bound $\mathds{1}_K$ in a different way than in Section \ref{sec:compare_grid_sdp}. We use an optimization problem in the form of \eqref{opt:boundSDPrbf}, without the constraints $X_i\succcurlyeq 0, \tau_i\geq 0$, thus getting an upper bound for $\mathds{1}_K$. This bound will not be as tight as the one computed with the gridding method of the previous example and as a result, all subsequent value function bounds will be less accurate. Our purpose here, however, is to compare controller performance and hence any bound for $\mathds{1}_K$ will work.

%\footnote{This may be a useful conclusion for a realistic system but does not yield useful results in the comparison we are making here.}

Despite having a different control objective, we assume that using an LQG controller where we choose the state and control penalty matrices according to the safe, target and control constraint sets is an efficient heuristic to maximize the reach-avoid probability. In particular we use the optimal control policy of the following problem:
\begin{align}
	\label{eq:lqr_problem_sdp}
		\begin{split}
	\min_{\bracks*{u_k}_{k=0}^{T-1}}\ &\mathbb{E}_{w_k}\parens*{\sum_{k=0}^{T-1}x_k^\top Q x_k + u_k^\top R u_k} + x_T^\top Q x_T
	\\
	&\st \quad\eqref{eq:linsys_sdp}, x_0 \in \bar{\Xs}
\end{split}
\end{align}
where $Q$ and $R$ have been chosen to correspond to the quadratic sets $K$ and $\Us$ respectively and penalize directions where they are small in size. Whenever the resulting LQG control input (calculated via the Riccati difference equation) is infeasible, we project it on the feasible set $\Us$. Starting from the same initial sates for both LQG and SDP-based controllers, we compare their performance by generating 100 trajectories using the same noise samples, counting the number that reach $K$ without leaving $K^\prime$. To make sure that the initial states chosen have a non-zero reach-avoid probability we first run the LQG controller and reject states with too low success probabilities. Note that since the design criteria of the two methods are different, our comparison is qualitative despite the fact that they both drive the system to the origin. Table \ref{tab:lqg_comp} reports the mean difference between the LQG and SDP-based controllers (denoted by $\hat{\pi}-\pi_{LQG}$) over the initial states for different problem dimensions. We notice that the SDP-based controller consistently outperforms the LQG controller in terms of success probability but scales worse in terms of computation time with the dimensions of the problem. The times reported for the LQG and SDP control, correspond to the computation of a single control input given the system state $x\in\Xs$. The last column of Table \ref{tab:lqg_comp} reports the time required to construct an upper bound value function for one step of the reach-avoid recursion using Algorithm \ref{algo:rbf_robust}. Note that a Gaussian RBF sum with a single element is needed to upper bound $\mathds{1}_K$ in each case ($M=1$). 

\begin{table}[t]
\centering
\begin{tabular}{ll|llll}
$n$&$m$&$\hat{\pi}-\pi_{LQG}$&LQG control (sec.) & SDP grid control (sec.) &SDP setup (sec.)\\
\hline
3 & 3 & 0.142 &  1.10e-04 & 0.0015&1.21 \\
3 & 4 & 0.126 &   9.951e-05 & 0.0302&1.11 \\
4 & 3 & 0.128 &   1.029e-04 &  0.0017&1.1\\
4 & 4 & 0.254 &   1.190e-04 & 0.0341 &1.09 \\
6 & 3 & 0.042 &    1.245e-04 & 0.0021 &1.08\\
5 & 4 & 0.163 &   1.215e-04 & 0.0430& 1.07 \\
6 & 4 & 0.072 &   1.174e-04 & 0.0504&1.094 \\
5 & 5 & 0.176 & 1.091e-04 & 0.846 & 1.07 \\
6& 5 & 0.144 & 1.256e-04  & 0.983 & 1.264 \\
7 & 5 & 0.133 & 1.43e-04 & 1.14 & 1.33  \\
8 & 5 & 0.141  & 1.33e-04 & 1.224 & 1.39\\
15 & 5 & 0.341  & 1.55e-04 &2.3 & 1.27 \\
\hline
\end{tabular}
\caption{Comparison of LQG and SDP-based controllers based on recursion \eqref{eq:DP} and Algorithm \ref{algo:rbf_robust}.}
\label{tab:lqg_comp}
\end{table}

\section{Conclusion}

The method presented in this paper is suitable for constructing upper bounds for the reach-avoid value function of Markov decision processes whose stochastic kernel can be expressed as a sum of Gaussian radial basis functions (RBFs) with affine dependence on the state and control input. In the case where safe and target sets can be written as intersections of quadratic inequalities, we derive conditions that satisfy the robust constraints of an infinite dimensional optimization problem,  equivalent to the reach-avoid dynamic programming recursion. We then provide a sequence of semidefinite programs to recursively construct upper bounds for the reach-avoid probability which are computationally efficient and scale to problems with much higher state and input dimensions than what can be handled in the literature. Our numerical examples indicate that the conditions used to reformulate the robust constraints may be conservative but the associated control policies perform close to optimal in low dimensional problems and considerably outperform a tuned LQG-based heuristic controller in higher dimensional problems.

As a part of future research we will investigate methods to reduce the conservatism introduced by the suggested constraint reformulation and couple the dominance conditions between Gaussian RBF sum functions. Moreover, it is worth comparing the suggested method to similar methods based on quadratic and polynomial approximations of value functions of general dynamic programs. We expect our method to provide a compromise between the scalability of quadratic approximations and the accuracy of general polynomial ones.

\begin{ack}
The authors would like to thank Marcello Colombino and Andreas Hempel from the Automatic Control Laboratory in ETH Zurich for discussions around the problem formulation. Special thanks to Mr. Colombino for his help in the reformulation of the constraints in the upper bound semidefinite program. 
\end{ack}

\bibliographystyle{asme_bibstyle}
\bibliography{ifasdp}

\end{document}